\documentclass[12pt]{article}
\usepackage{amssymb}

\def\appendix{\par}  
\def\poset{{\mathbb P}}
\def\uu{{\mathcal U}}
\def\vv{{\mathcal V}}
\def\ww{{\mathcal W}}
\def\hh{{\mathbb H}}
\def\ll{{\mathbb L}}
\def\cc{{\mathcal C}}
\def\bb{{\mathfrak b}}
\def\split{{\rm split}}
\def\root{{\rm root}}
\def\leaves{{\rm leaves}}
\def\splitnode{{\rm splitnode}}
\def\diamondb{\Diamond(\bb)}
\def\om{\omega}
\def\infsub{[\om]^\om}
\def\finsub{[\om]^{<\om}}
\def\sq{\subseteq}
\def\lq{\mbox{``}}
\def\rq{\mbox{''}}
\def\proof{\par\noindent Proof\par\noindent}
\def\name#1{\stackrel{\circ}{#1}}
\def\reals{{\mathbb R}}
\def\forces{{| \kern -2pt \vdash}}
\def\force{\forces}
\def\res{\upharpoonright}
\def\qed{\par\noindent QED\par}
\def\rmand{\mbox{ and }}
\def\rmor{\mbox{ or }}

\newtheorem{theorem}{Theorem}
\newtheorem{lemma}[theorem]{Lemma}
\newtheorem{cor}[theorem]{Corollary}

\newtheorem{question}[theorem]{Question}

\newtheorem{prop}[theorem]{Proposition}

\begin{document}

\begin{center}
{\large The $\gamma$-Borel conjecture}
\end{center}

\begin{flushright}
Arnold W. Miller\footnote{
Thanks to Boise State University  
for support during the time this paper
was written and to Alan Dow for some helpful discussions and to
Boaz Tsaban for some suggestions to improve an earlier version.  
\par Mathematics Subject Classification 2000: 03E35; 03E17
\par Keywords: $\gamma$-set, Hechler forcing,
 Laver forcing, strong measure zero set.
}
\end{flushright}

\begin{center}
Abstract
\end{center}

\begin{quote}
In this paper we prove that it is consistent that every $\gamma$-set
is countable while not every strong measure zero set is countable.
We also show that it is consistent that every strong $\gamma$-set is
countable while not every $\gamma$-set is countable.  On the other
hand we show that every strong measure zero set is countable iff
every set with the Rothberger property is countable.
\end{quote}

A set of reals $X$ has strong
measure zero iff for any sequence $(\epsilon_n:n<\om)$ of
positive reals there exists a sequence of intervals $(I_n:n<\om)$
covering $X$ with each $I_n$ of length less
than $\epsilon_n$.  Laver \cite{laver} showed that it is relatively
consistent with ZFC that the Borel conjecture is true, i.e.
every strong measure zero set is countable. 

Sets of reals called $\gamma$-sets were first considered by Gerlits  and Nagy
\cite{gn}.  They showed that every $\gamma$ set has strong measure zero and
that Martin's Axiom implies every set of reals of size smaller than the
continuum is a $\gamma$-set.   A $\gamma$-set of size continuum is constructed
in Galvin and Miller \cite{gm} using MA.

Next we define $\gamma$-set.  An open cover  $\uu$ of a topological space $X$
is an $\om$-cover iff  for every finite $F\sq X$ there exists $U\in \uu$ with
$F\sq U$ and $X\notin\uu$.  An open cover $\uu$ of $X$ is a $\gamma$-cover iff
$\uu$ is infinite and each $x\in X$ is in all but finitely many $U\in\uu$. 
Finally, $X$ is a $\gamma$-set iff $X$ is a separable metric space in which
every $\om$-cover contains a $\gamma$-subcover.

\bigskip

Paul Szeptycki asked if it was possible to have a sort of weak Borel
conjecture  be true, i.e., every $\gamma$-set countable, while the Borel
conjecture is false.  We answer his question in the positive.  We use  Hechler
\cite{hechler} forcing, $\hh$, for adding a dominating real, an analysis of it
due to Baumgartner and Dordal \cite{bd}, and properties of
Laver forcing $\ll$, and a characterization
of $\hh$ due to Truss \cite{truss}.

\bigskip

\begin{theorem} \label{mainthm}
If $\hh$ is iterated $\om_2$ times with finite support
over a model of CH, then in the resulting model every $\gamma$-set
is countable but every set of reals of cardinality $\om_1$ has
strong measure zero.
\end{theorem} 

\proof

For $f\in \om^\om$, define $\uu_f$ to be the following family of
clopen subsets of $2^\om$. 
 $$\uu_f=\{C_F:\exists n\;\; F \sq 2^{f(n)}, |F|\leq n\} \mbox{ where }
 C_F=\{x\in 2^\om: x\res n\in F\}\}$$
Note that for any finite $A\sq 2^\om$ there exists $C\in \uu_f$ with
$A\sq C$.  Also $2^\om\notin U_f$ provided that $f(n)>2^n$ all $n$.
Let $\ll$ denote Laver forcing \cite{laver}.

\begin{lemma}\label{one} Suppose $M$ is a model of set theory, $f$ is 
$\ll$-generic over $M$, and $X\sq 2^\om$ is in $M$.  Then 
$$M[f]\models\forall \cc\in [\uu_f]^\om\;\;
|\bigcap{\cc}\cap X |\leq \om $$ 
\end{lemma}
\proof

For a tree $p\sq \om^{<\om}$ and $s\in p$ we define
$$p_s=\{t\in p: t\sq s \rmor s\sq t\}$$

A Laver condition (or Laver tree) is a tree $p\sq \om^{<\om}$
with a root $s\in p$ with the property that $p_s=p$ and
for every $t\in p$ with $|t|\geq |s|$ there exists infinitely many
$n<\om$ with $tn\in p$. The order is $p\leq q$ iff $p\sq q$.
As usual we define $p\leq_0 q$ iff $p\leq q$ and $\root(p)=\root(q)$.
Somewhat nonstandardly let us write 
$$\leaves(p)=\{r\in p: \root(p)\sq r\}$$
and for each $s\in\leaves(p)$ define
$$\split(p,s)=\{n\in\om: sn\in p\} $$

Suppose that the lemma is false. 
Let $p$ be a Laver condition such that  
$$p\forces\lq 
\bigcap\name{\cc}\cap X =\name{Y}\mbox{ is uncountable }
\rmand \name{\cc}\sq \uu_f \mbox{ is infinite}\rq$$

By cutting $\cc$ down (if necessary)  we may suppose that 
$\cc=\{C_{F_n}: n\in Q\}$ where $F_n\sq 2^{f(n)}$ with  $|F_n|\leq n$
and  $Q\in \infsub$.

Working in $M$ using standard arguments
of Laver forcing \cite{laver} we can prove the following 
Claims.

\bigskip
{\bf Claim.} Suppose that $p$ is an arbitrary condition such that 
$$p\forces \{\name{s}_i:i<k\}\sq 2^{f(k)}$$ 
where $k=|s|$ and $s=\root(p)$. Then there exists $r\leq_0 p$ and   
$(x_i\in 2^\om:i<k)$ such that for any $m<\om$ 
for all but finitely
many $n\in\split(r,s)$ for every $i<k$
$$r_{sn}\forces x_i\res m = s_i\res m$$ 

\proof 
One of the basic properties of Laver forcing is that if $p$ is any Laver tree
and $\theta$ any sentence in the forcing language, then there exists 
$q\leq_0 p$, which decides $\theta$, i.e.     
$$q\forces \theta \;\; \rmor \;\;q\forces\neg\theta.$$   

Note that for $sn\in p$ we
have that $p_{sn}\forces f(k)=n$.  Hence we can find $q\leq_0 p$ and
$(s_i^n\in 2^n : i<k, n\in\split(q,s))$ so that for each $n\in\split(q,s)$ 
we have that

$$q_{sn}\forces \lq s_i^n=\name{s}_i \mbox { for all } i<k\rq$$ 

It follows by compactness that there exists $x_i\in 2^\om$
and an infinite set $E\sq\split(q,s)$ so that for every
$m<\om$ we have that for all but finitely many $n\in E$ that
 $$s_i^n\res m = x_i\res m \mbox{ for all } i<k .$$  
Now let $r=\cup \{q_{sn}:n\in E\}$ so that $r\leq_0 q $. 
\qed

Note that if
$y\in 2^\om\setminus\{x_i:i<k\}$, then 
$$r_{sn}\forces y\notin C_{\{s_i:i<k\}}$$
for all but finitely
many $n\in \split(r,s)$. 
By the usual fusion arguments we obtain:

\bigskip
{\bf Claim} There exists  $q\leq_0 p$ and 
$(K_s\in [2^\om]^{\leq |s|}:s\in \leaves(q))$ such that
\begin{enumerate}
 \item for each $s\in \leaves(q)$ either
  $q_s\forces |s|\in \name{Q}$ or $q_s\forces |s|\notin \name{Q}$, and
 \item for each $s\in\leaves(q)$ if $q_s\forces |s|\in \name{Q}$
 then for any $x\in 2^\om\setminus K_{s}$ for all but finitely
 many $n$ if $sn\in q$, then 
 $$q_{sn}\forces x\notin  \name{C}_{{F}_{|s|}}$$
\end{enumerate} 

\proof
We repeat the first Claim at each node starting at the root and continuing
downward and then take the fusion.
\qed

\bigskip

Now since $p$ forces that ${Y}$ is uncountable we must be able to find 
$$x\notin\cup \{K_s\;:\; s\in \leaves(q)\}$$ 
and $r\leq q$ such that 
$r\forces x\in Y$.  But this is a contradiction, since there must
be some $s\in\leaves(r)$ such that 
$$r_s\forces\lq |s|\in \name{Q}\rq$$
 and then for all but finitely many 
$n\in\split(r,s)$ we have that 
$$r_{sn}\forces x\notin\name{C}_{{F}_{|s|}}$$
But even one such $n$ gives a contradiction. This proves the Lemma.
\qed

Now we note that this property is preserved when we add a Cohen real.

\begin{lemma} \label{two} Suppose $N$ is a model of set theory, 
$x\in\om^\om$ is a Cohen real
over $N$, $X\sq 2^\om$ in  $N$  and 
$\uu\in N$ is a family of subsets of $X$, and 
$$N\models\lq \forall \cc\in [\uu]^\om\;\;
|\bigcap{\cc}\cap X|\leq \om \rq$$ 
Then 
$$N[x]\models\lq \forall \cc\in [\uu]^\om\;\;
|\bigcap{\cc}\cap X|\leq \om$$
\end{lemma}

\proof

Suppose not and let 
$$p\force X\cap\bigcap\name{\cc} \mbox{ is uncountable }$$
Since the Cohen partial order is countable, there would
exist $q\leq p$ so that 
$$Y=\{x\in X: q\forces x\in \bigcap\name{\cc}\}$$
is uncountable and in $N$.  But then letting
$$\cc'=\{U\in\uu: Y\sq U\}$$
yields a contradiction.
\qed

It follows from the two Lemmas that if $f$ is Laver over $M$, 
$x$ is Cohen over $N=M[f]$, and $X\sq 2^\om$ is an uncountable
set in $M$, then in $M[f,x]$ every infinite $\cc\sq\uu_f$ 
has the property that $\bigcap\cc\cap X$ countable.

The following Lemma applies to the Laver real $f$ since it is
dominating.

\begin{lemma} (Truss \cite{truss}) Suppose $f$ is a dominating real over $M$,
i.e., $g\leq^*f$ for every $g\in M\cap\om^\om$ and 
$x\in \om^\om$ is a Cohen real over $M[f]$, then 
$h=f+x$ is $\hh$-generic over $M$.  
\end{lemma}

\begin{lemma} \label{main} Let $f$ be $\ll$-generic over $M$, \label{five}
$x\in\om^\om$ a Cohen real over $M[f]$, and $h=f+x$. 
Then for every uncountable
$X\sq 2^\om$ in $M$ 
$$M[h]\models\lq \forall \cc\in [\uu_h]^\om\;\;
\bigcap{\cc}\cap X \mbox{ is countable }\rq.$$ 
\end{lemma}
\proof
Note that in $M[f,x]$ that every infinite  $\cc\sq\uu_h$ 
has the property that $\bigcup\cc\cap U_f$ countable.  To see this
suppose otherwise and consider $\{C_{H_n}:n\in Q\}$ with 
$H_n\sq 2^{h(n)}$, $Q$ infinite, and the
$C_{H_n}$ distinct.  Define $F_n=\{s\res f(n): s\in H_n\}$.
Now since $h(n)\geq f(n)$ we have that $C_{H_n}\sq C_{F_n}$.
The set $\{C_{F_n}:n\in Q\}$ must be infinite because 
the measure of $C_{F_n}$ is $\leq {n\over{2^{f(n)}}}$.   
Since $M[h]\sq M[f,x]$ the lemma follows.
\qed  

Note that the lemma applies to every Hechler generic real and not just the sum
of a Laver and a following Cohen.  This is because if it is false it must be
forced false by a particular Hechler condition.  
Then just take a Laver real in
that condition and follow it with a Cohen to get a contradiction.
In more detail let 
$$\hh=\{(n,f):f\in\om^\om, n\in\om\}$$
and define the Hechler neighborhoods
$$[n,f]=\{g\in\om^\om: g\res n= f\res \rmand \forall i \;\;g(i)\geq f(i)\}$$
Then $(m,g)\leq (n,f)$ iff $m\geq n$ and $g\in [n,f]$.  Also
for $G$ $\hh$-generic over $M$ the Hechler real is 
$$h=\bigcup\{f\res n:(n,f)\in G\}$$
and it has the property that
$$G=\bigcup\{(n,f)\in \hh: h\in [n,f]\}$$
The lemma must be true in every Hechler extension,
If not, there would exist some condition
$(n,g)$ forcing it is false.  It is easy to find a Laver real 
$f\in [n,g]$ and letting $x\in\om^\om$ be a Cohen real over
$M[f]$ with $x\res n$ constantly zero, we would get a Hechler real
$h=f+x$ with  $h\in [n,g]$ which gives a contradiction.

\begin{question}(Ramiro de la Vega) Given a countable transitive model
of set theory $M$, is it true that for every Hechler real $h$
over $M$ there exists a Laver real $f$ over $M$ and a Cohen real
$x$ over $M[f]$ such that $h=f+x$?
\end{question}
 
\bigskip 
Define $(a_\alpha\in \infsub:\alpha<\om_1)$ is eventually narrow iff
for every $b\in\infsub$ there exists $\alpha<\om_1$ so that $b\setminus
a_\beta$ is infinite for all $\beta>\alpha$.

\begin{lemma} (Baumgartner and Dordal \cite{bd}) Suppose $N$ is
a model of set theory and 
$$N\models (a_\alpha\in \infsub:\alpha<\om_1)
\mbox{ is eventually narrow. }$$
Then for any $G_{\om_2}$ which is $\hh_{\om_2}$-generic
over $N$, we have that 
$$N[G_{\om_2}]\models (a_\alpha\in \infsub:\alpha<\om_1)
\mbox{ is eventually narrow. }$$
\end{lemma}  

Now we prove that every $\gamma$-set in $M[G_{\om_2}]$ countable. 
Since $\gamma$-sets are zero dimensional we need only worry about
uncountable $Y\sq 2^\om$. Let $X\sq Y$ be a subset
of size $\om_1$.  Construct $g:\om\to\om$ so that for every
$n<\om$ if $m=g(n)$, then 
$$|\{x\res m:x\in X\}|>n$$
By the usual ccc finite support iteration arguments we can find
$\alpha<\om_2$ so that $X,g\in M[G_\alpha]$ and  letting $h=h_\alpha$ be
the next Hechler real added we
have that $h(n)>g(n)$ for all $n$. 
From Lemma \ref{main} and the remark following it we that in 
$N=M[G_{\alpha+1}]$
for every infinite $\cc\sq \uu_h$ that
$\bigcap\cc\cap X$ is countable.  
Now since $h(n)>g(n)$ there is no $U\in U_h$ which covers $X$, however
$\uu_h$ is an $\om$-cover of $2^\om$ and hence of $Y$.

Now let $X=\{x_\alpha:\alpha<\om_1\}$ and $\uu_h=\{U_n:n<\om\}$.
In the model $N=M[G_{\alpha+1}]$
define $a_\alpha=\{n<\om: x_\alpha\in U_n\}$.  Note that
$$N\models (a_\alpha\in \infsub:\alpha<\om_1)
\mbox{ is eventually narrow. }$$

Otherwise if $b\sq^* a_\alpha$ for uncountably many  $\alpha$, then 
for some infinite $c\sq b$ 
$$Y=\{x_\alpha: c\sq a_\alpha\}$$
is uncountable. But then $Y\sq \bigcap\{U_n:n\in c\}$ 
which  contradicts Lemma \ref{five}.

Since the tail of a finite iteration of $\hh$ is itself a finite
support iteration $\hh$ the
Baumgartner-Dordal Lemma applies and so, 
 $$N[G_{[\alpha+2,\om_2)}]= M[G_{\om_2}]$$ 
models that $(a_\alpha:\alpha<\om_1)$ is eventually 
narrow.
But this implies that $Y$ is not a $\gamma$-set since if 
$(U_n\in\uu:n\in b)$ is a $\gamma$-cover of $X\sq Y$, then 
for some infinite $c\sq b$, we would have that
$X\cap \bigcap\{U_n:n\in c\}$ is uncountable, which implies that
for uncountably many $\alpha$ that $c\sq a_\alpha$. Contradicting
the fact the $a_\alpha$ are eventually narrow.

On the other hand, it is well known that forcing with $\hh$ adds
Cohen reals and adding Cohen reals makes
sets of reals of small cardinality into strong measure zero sets.  
To see this suppose that
$(\epsilon_n>0:n<\om)\in M$ a model of set theory.  In $M$ let
$(I_{nm}:m<\om)$ list all intervals with rational end points and
of length less than $\epsilon_n$. If
$x:\om\to \om$ is a Cohen real over $M$, then it is an easy density
argument to prove that
$$M\cap\reals \sq \bigcup_{n<\om} I_{n x(n)}$$

The usual arguments show that in the iteration every set of reals of
cardinality $\om_1$ has strong measure zero.
This proves Theorem \ref{mainthm}.
\qed

\bigskip Remark. It is also true in the Hechler real model that every set of
reals of size $\om_1$ is both in ${\sf S}_1(\Gamma,\Gamma)$ and  
${\sf S}_1(\Omega,\Omega)$.   For definitions, see Just, Miller,  Scheepers, 
and Szeptycki  \cite{jmss}.  This follows from the fact that $\bb>\om_1$ 
and $cov(\mathcal{M})>\om_1$, see Figure 4 \cite{jmss}.

\bigskip
Define. $X$ is $C''$  iff for every sequence
$(\uu_n:n<\om)$ of open covers of $X$ there exist 
$(U_n\in\uu_n:n<\om)$ an open cover of $X$.
Equivalent terminology for $C''$ is the Rothberger property or 
${\sf S}_1({\mathcal O},{\mathcal O})$.

Define. $C''$-BC to be the statement that every set of reals with the
property $C''$ is countable and let SMZ-BC denote the standard Borel
conjecture, every strong measure zero set is countable.

\begin{prop}
SMZ-BC is equivalent to $C''$-BC.
\end{prop}
\proof
It is only necessary to prove right to left.

If $\bb=\om_1$ then there exists an uncountable set of reals concentrated
on the rationals (Rothberger) and any such set has property $C''$.
So assume $\bb>\om_1$.

Suppose there is an uncountable strong measure zero set.  Then by
standard arguments there exists an
$X\sq 2^\om$ with $|X|=\om_1$ such  that for every $f\in\om^\om$ there exists
$(s_n\in 2^{f(n)}:n<\om)$ such that for every $x\in X$ there are infinitely
many $n$ with $s_n\sq x$.

\bigskip\noindent
{\bf Claim.} $X$ has property $C''$. 

\proof
Let $(\uu_n:n<\om)$ be open covers of $X$.  Without loss we may assume
each element of each $\uu_n$ is of the form $[s]$ for some $s\in 2^{<\om}$.
Since $|X|<\bb$ we can find finite $A_n\sq 2^{<\om}$ so that
$s\in A_n$ implies $[s]\in \uu_n$ and for each $x\in X$ for all but
finitely many $n$ there exists $s\in A_n$ with $s\sq x$.  Let
$f:\om\to\om$ be such that $f(n)>\max\{|s|:s\in A_n\}$.  Using
strong measure zero of $X$ choose $s_n\in 2^{f(n)}$ so that
every element of $X$ is in infinitely many $[s_n]$.  Define $t_n\in A_n$
as follows. If there exists $t\in A_n$ with $t\sq s_n$ then
let $t_n$ be such.  If there isn't, choose $t_n$ arbitrarily.  We
claim that $\{[t_n]:n<\om\}$ covers $X$. For any $x\in X$ for all
but finitely many $n$ we have that there exists $t\in A_n$ with
$t\sq x$.  But for infinitely many $n$ we have that $s_n\sq x$.
Since $|s_n|>|t_n|$ it must be the case that for infinitely many $t_n$ that
$t_n\sq x$.

This proves the Claim and the Proposition.
\qed 

\bigskip 

Define $X$ is a strong $\gamma$-set iff there exists an increasing sequence of
integers $(k_n:n<\om)$ so that for every sequence $(\uu_n:n<\om)$ where $\uu_n$
is a $k_n$-cover of $X$ (i.e. covers every $k_n$ element subset of $X$) there
exists a $\gamma$-cover of the form $(U_n\in\uu_n:n<\om)$.   These were first
defined in Galvin and Miller \cite{gm}.  Tsaban \cite{tsaban} has shown
that an equivalent definition results if we always require $k_n=n$.

\begin{theorem}
In the Cohen real model, i.e., $\om_2$ Cohen reals added to a model
of CH, every strong $\gamma$-set is countable but there is an uncountable
$\gamma$-set.
\end{theorem}
\proof

First we construct an uncountable $\gamma$-set.
This proof is a modification of the construction ??? from Just, Miller,
Scheepers, and Szeptycki \cite{jmss}.

Without loss of generality we may assume that 
$$N=M[x_\alpha\sq\om:\alpha<\om_1]$$ 
where  the Cohen reals occur at the end.  Note $M$ fails to
satisfy CH. Construct  $y_\alpha\in\infsub$ descending mod finite so that
$(y_\beta: \beta<\alpha)\in M[x_\beta:\beta<\alpha]$ as follows:  

At stage $\alpha+1$ let 
$$y_{\alpha+1}=x_{\alpha+1}\cap y_\alpha$$ 
This is infinite because $x_{\alpha+1}$ is Cohen generic over $y_\alpha$ At
limit stages choose $y_\alpha\in M[x_\beta:\beta<\alpha]$ in
some canonical way (maybe using sequence of enumeration of the
countable ordinals in $M$) so that
$y_\alpha\sq^*y_\beta$ all $\beta<\alpha$.

\bigskip\noindent
{\bf Claim.} Suppose $(\uu_n:n<\om)\in M[x_\beta:\beta\leq \alpha]$ 
is a family of
$\om$-covers of 

$$\finsub\cup \{y_\beta:\beta\leq \alpha\}$$
 
Then there exists a sequence $(U_n\in\uu_n:n<\om)$
which is a $\gamma$-cover of  

$$\finsub\cup \{y_\beta:\beta\leq \alpha\} \cup
[y_{\beta+1}]^{*\om}$$

\bigskip
\proof

Let $\bigcup_n F_n=\finsub\cup \{y_\beta:\beta\leq \alpha\}$
be an increasing union of finite sets and define
$\vv_n=\{U\in\uu_n : F_n\sq U\}$
and note that they are $\om$-covers.
Next inductively define $\ww_n$ by $\ww_0=\vv_0$ and 
$$\ww_{n+1}=\{U\cap V: U\in\vv_n, V\in\ww_n\}$$ 
and note that they are $\om$-covers which refine each other. 
Working in the ground model construct an increasing sequence 
$k_n$ and $U_n\in \ww_n$ so that 
$$\{x\sq\om:x\cap [k_n,k_{n+1})=\emptyset\}\sq U_n$$
this can be done since $\ww_n$ is an $\om$-cover
of $\finsub$.  Now since $x_{\alpha+1}$ is Cohen
real the following set will be infinite:
$$A=\{n<\om: x_{\alpha+1}\cap [k_n,k_{n+1})=\emptyset\}$$ 
The same or larger set will work for $y_{\alpha+1}$
and so $(U_n:n\in A)$ will be a $\gamma$-cover of
$[y_{\beta+1}]^{*\om}$. The refining conditions on $\ww_n$ means we can
fill it in on the complement of $A$ and the choice of $\vv_n$ means
it is a $\gamma$-cover of the rest.  
\qed

The Claim shows that
$\finsub\cup\{y_\alpha:\alpha<\om_1\}$ is
a $\gamma$-set.

\bigskip Next we show that there are no uncountable strong $\gamma$-sets.
Suppose for contradiction that $X\sq 2^\om$ is an uncountable strong
$\gamma$-set witnessed by $(k_n:n<\om)$ in the model $N$. By the usual ccc
arguments we may suppose that $X,(k_n:n<\om)\in M$ where $M\sq N$ is some model
of CH.  Let $u\in N\cap \om^\om$ be Cohen generic over $M$ and $v\in N\cap
\om^\om$ Cohen generic over $M[u]$ so that if we let 
$$\uu_n=\{[s]:s\in 2^{u(n)}\}$$
then (using that $N$ thinks $X$ is strong $\gamma$) there exists 
 $$(\vv_n\in [\uu_n]^{\leq k_n}:n<\om)\in M[u,v]$$
so that 
$\forall x\in X\forall^\infty n\;\; x\in\cup\vv_n$.  
Let $\poset$ denote Cohen forcing and since
it is countable there must be some $(p,q)\in\poset\times\poset$ and 
$N<\om$ such that 
$$(p,q)\forces (\name{\vv}_n\in [\name{\uu}_n]^{k_n}:n<\om)$$
and
$$Y=\{x\in X:(p,q)\forces \forall n>N \;\; x\in\cup\name{\vv}_n\}$$
is uncountable.  Fix $n>N,|p|$.  Now since $Y$ is uncountable there
exist some level $l<\om$ with 
$$|\{x\res l: x\in Y\}|> k_{n}$$
Let $r\supseteq p$ be an extension with $r(n)=l$.  But this is a
contradiction since
\begin{itemize}
 
 \item $(r,q)\forces \lq\vv_n\sq 2^{l} \rmand |\vv_n|\leq k_{n}\rq$, and
 
 \item $(r,q)\forces \lq x\in\vv_{n}\rq$ for every $x\in Y$ and
 so $(r,q)\forces \lq \{x\res l: x\in Y\}\sq \vv_n\rq$
  
\end{itemize}
\qed
  
Remark. T. Bartoszynski has shown that in the iterated superperfect real model
every strong $\gamma$-set is countable.  Superperfect forcing  is also called
rational perfect set forcing, see Miller \cite{super}.   The principle
$\diamondb$ (see Dzamonja, Hrusak, and Moore \cite{dhm}) implies that there is 
an uncountable $\gamma$-set.  Since $\diamondb$ holds in the iterated
superperfect real model, we get another model for the consistency of strong
$\gamma$-BC but not $\gamma$-BC.

\begin{flushleft}
Arnold W. Miller \\
miller@math.wisc.edu \\
http://www.math.wisc.edu/$\sim$miller\\
University of Wisconsin-Madison \\
Department of Mathematics, Van Vleck Hall \\
480 Lincoln Drive \\
Madison, Wisconsin 53706-1388 \\
\end{flushleft}

\appendix

\newpage
\begin{center}
Appendix \\ 
\end{center}

This is not intended for publication but only for the electronic version.

\begin{theorem}(T. Bartoszynski)
In the iterated superperfect forcing model, every strong $\gamma$-set
is countable. 
\end{theorem}

\proof
This model is obtained by the countable support iteration of length
$\om_2$ of superperfect forcing  over a model of
CH.

First we consider one-step. Let
$f$ be superperfect generic over $M$ a model of set theory. Define 
$(\uu_n:n<\om)$ by $$\uu_n=\{[s]: s\in 2^{f(n)}\}.$$

\bigskip\noindent{\bf Claim}. Let $g\in \om^\om\cap M$ and
$(\vv_n\in [\uu_n]^{<g(n)}:n<\om)\in M[f]$.
Then 
$$M[f]\models 
|\{x\in M\cap 2^\om\;:\; \forall^\infty n\;\; x\in \cup\vv_n)|\leq\om.$$

\proof
For $p$ a superperfect tree, define $s\in\splitnode(p)$ iff
$\exists^\infty n\; sn\in p$.  Superperfect trees are those trees
in which the split nodes are dense.  Suppose 
$$p\forces (\name{\vv}_n\in [\uu_n]^{<g(n)}:n<\om)$$

By the usual fusion arguments we can obtain a superperfect tree $q\leq p$ and 
$(K_s\sq 2^\om\;:\; s\in\splitnode(q))$ so that
\begin{enumerate}
 \item $|K_s|< g(|s|)$ for each $s\in\splitnode(q)$ 
 \item for each $s\in\splitnode(q)$ and $x\in 2^\om\setminus K_s$ 
 for all but finitely many $n\in\split(q,s)$  
 $$q_{sn}\forces x\notin\name{\vv}_{|s|}$$
\end{enumerate}
It follows that
$$ q\forces \lq M\cap(\cup_{m<\om}\cap_{n>m}\cup \name{\vv}_n)\sq
\cup\{K_s\;:\; s\in\splitnode(q)\}\rq$$

\qed
Now suppose for contradiction that $X$ is an uncountable 
strong $\gamma$-set in the model $M[f_\alpha:\alpha<\om_2]$.
By the $\om_2$ chain condition and a Lowenheim-Skolem argument
there must be an $\alpha_0<\om_2$ 

with $X,(k_n:n<\om)\in M[f_\alpha:\alpha<\alpha_0]$ such that
$$M[f_\alpha:\alpha<\alpha_0]\models X \mbox{ is a strong $\gamma$-set
with witness } (k_n:n<\om)$$

Denote $M[f_\alpha:\alpha<\alpha_0]$ as $M_0$.
Now using $f_{\alpha_0}$ (the next superperfect real) 
Let $\uu_n=\{[s]\;:\; s\in 2^{f_{\alpha_0}(n)}\}$. By the one step argument
for any $g\in M_0\cap \om^\om$ 
$$M_0[f_{\alpha_0}]\models\forall (\vv_n\in[\uu_n]^{g(n)}
|\{x\in X\cap 2^\om\;:\; \forall^\infty n\;\; x\in \cup\vv_n)|\leq\om.
$$

Denote $M_0[f_{\alpha_0}]$ as $M_1$.  Our final model 
$M_2 = M[f_\alpha:\alpha<\om_2]$ satisfies the Laver property over
the intermediate models. .  This means 
for any $f\in M_2\cap \om^\om$ such that there exists
$h\in M_1\cap \om^\om$ which bounds $f$, i.e., $f(n)<h(n)$ all $n$, there
exists $(H_n:n<\om)\in M_1$ with $|H_n|\leq 2^n$ and $f(n)\in H_n$ for
all $n$.   The reason this is true is that the Laver property holds
in the one-step superperfect model by essentially the same argument as
for Laver forcing. It also holds in the iteration by either the same
argument Laver employed or by the general fact that it is preserved
by countable support iteration of proper forcings (see Bartoszynski
and Judah; {\bf Set theory. On the structure of the
real line.} A K Peters, Ltd., Wellesley, MA, 1995.).

But now we get a contradiction.
Let $\uu_n^*$ to be the family of $k_n$ unions of elements of $\uu_n$. 
Since $M_2$ thinks that $X$ is a strong $\gamma$-set there
is a $\gamma$-cover of $X$ of the form $(V_n\in\uu_n^*:n<\om)$.
But by the Laver property this means there exists
$(\vv_n\in [\uu_n]^{k_n 2^n}:n<\om)\in M_1$ with $V_n\sq \cup \vv_n$.
But this is a contradiction for $g(n)=k_n 2^n$ and $M_1=M_0[f_{\alpha_0}]$.  

\qed
Next we show that there is an uncountable $\gamma$-set in the superperfect
model.  We construct it using the principle $\diamondb$.  This is stronger than
$\bb=\om_1$ and is defined in Dzamonja, Hrusak, and Moore \cite{dhm}.   They
prove that it holds in any model of  $\bb=\om_1$ which is obtained by the
$\om_2$-iteration with countable support of proper Borel orders which
are reasonably homogeneous. Hence, $\diamondb$ is true in the iterated
superperfect set forcing model.

I do not know if $\bb=\om_1$ is enough to construct an uncountable 
$\gamma$-set.

\bigskip 
Define $\diamondb$: For every $F:2^{<\om_1}\to \om^\om$ 
such that each $F\res 2^\alpha$ is Borel for $\alpha<\om_1$
there exists $g:\om_1\to
\om^\om$ so that for every $f\in 2^{\om_1}$ such 
$\exists^{\infty} n \; F(f\res\delta)(n)<g(\delta)(n)$ for stationarily
many $\delta<\om_1$. 
\bigskip

\begin{theorem}
$\diamondb$ implies there is an uncountable $\gamma$-set.
\end{theorem}

\proof

Let  $H_\delta:(\infsub)^\delta\to \infsub)$ be Borel so that
for any $(x_\alpha:\alpha<\delta)$ if 
and $\alpha<\beta$ implies $x_\beta\sq^* x_\alpha$, then
for $y=H(x_\alpha:\alpha<\delta)$ we have that $y\sq^*x_\alpha$
for every $\alpha<\delta$.
By using the first $\om$-coordinates to code a countable family
of open sets we may assume that the domain of $F$ is 
sets of the form 
$(\uu_n:n<\om), (x_\alpha\sq\om:\alpha<\delta)$ where
the $\uu_n$ are families of open subsets of $2^\om$ and
we are to define 
$$F((\uu_n:n<\om), (x_\alpha:\alpha<\delta))=h\in\om^\om$$

Suppose 
\begin{enumerate}
\item $x_\alpha\in\infsub$ for each $\alpha<\delta$,
\item $x_\alpha\sq^* x_\beta$ for each $\beta<\alpha<\delta$, and
\item $\uu_n$ is an $\om$-cover of $\finsub\cup\{x_\alpha:\alpha<\delta\}$
for each $n$. 
\end{enumerate} 
(If any of these fail to be true, just define $h$ to be the constant zero
function.) 

Let $\{\delta_i:i<\om\}=\delta$ 
be some previously chosen 
enumeration of $\delta$ and define for each $n$
$$\vv_n=\{U\in\uu_n:\{x_{\delta_i}:i<n\}\sq U\}$$
It is easy to check that each $\vv_n$ is an $\om$-cover of 
$\finsub\cup\{x_\alpha:\alpha<\delta\}$.  Also choosing an element of
each will automatically $\gamma$-cover $\{x_\alpha:\alpha<\delta\}$.
Next define inductively $\ww_n$ as follows:
\begin{enumerate}
\item $\ww_0=\vv_0$,
\item $\ww_{n+1}=\{U\cap V:U\in \ww_n, V\in \vv_n\}$
\end{enumerate}
It is easy to check that the intersections of elements of two $\om$-covers
is an $\om$-cover, so by induction each $\ww_n$ is an $\om$-cover
of $\finsub\cup \{x_\alpha:\alpha<\delta\}$.  Since $\ww_{n+1}$ is
a refinement of $\ww_n$, if for some $A\in \infsub$ we have
$(U_n\in \ww_n:n\in A)$ is a $\gamma$-cover, then we can choose $U_n$
for $n\notin A$ by looking forward to the next element of $A$ so that
$(U_n\in \ww_n:n\in\om)$ is a $\gamma$-cover.  

Apply $H$ to get $H(x_\alpha:\alpha<\delta)=\{k_n:n<\om\}$ (Note
that this does not depend on the covers $\uu_n$.)  Construct an
infinite $B\sq \om$ so that for every successive pair of elements
of $B$, say $n<m$, there exists $U_n\in \ww_n$ so that
$$\{x\sq\om: x\cap [k_n,k_m)=\emptyset\}\sq U_n$$
This only uses that $\ww_n$ is an $\om$-cover of $\finsub$: choose
$U$ to cover $[k_n]^{<\om}$ and then using that $U$ is open make
sure that $k_m$ is sufficiently large.  Now we make sure that 
$h\in\om^\om$ is such that $h$ eventually dominates the enumeration
function of $B\setminus N$ for each $N<\om$.  We leave to the reader
the details of showing that $h$ can be obtained using a Borel function on
$(\uu_n:n<\om), (x_\alpha:\alpha<\delta)$.  But note the following:
Suppose $g\in\om^\om$ has the property that $\exists^\infty n\;
g(n)>h(n)$, then there must be infinitely many $i$ so that there exists
$n<m$ elements of $B$ so that $g(i)\leq n<m \leq g(i+1)$.  Otherwise the
enumeration function of some $B\setminus N$ would dominate $g$ which is
impossible.

Applying $\diamondb$ to our function $F$ 
we get a $g:\om\to\om^\om$. Construct our $\gamma$-set
$X=\finsub\cup\{x_\alpha:\alpha<\om_1\}$ as follows:

Given $\{x_\alpha:\alpha<\delta\}$ and descending sequence in $\sq^*$
apply  $H$ to get $H(x_\alpha:\alpha<\delta)=\{k_n:n<\om\}$. Let
$g=g(\delta)\in\om^\om$ and put $x_\delta=\{k_{g(n)}:n<\om\}$.  
Now we verify that $X$ is a $\gamma$-set.  Suppose that $(\uu_n:n<\om)$
are open $\om$-covers of $X$.  By the definition of $\diamondb$ there
are stationarily many $\delta<\om_1$ such that
$$F((\uu_n:n<\om), (x_\alpha:\alpha<\delta))=h\in\om^\om$$
and if $g=g(\delta)$, then $\exists^\infty n\; g(n)>h(n)$.  
Note that $x_\delta=\{k_{g(i)}:n<\om\}$.  So as we
have remarked there are infinitely many $i$ (say $i\in C$) 
so that there exists
elements of $n_i<m_i$ of $B$ with $g(i)<n_i<m_i<g(i+1)$. The way
the elements of $B$ were construct means that there exists 
$U_{n_i}\in \ww_{n_i}$ such that
$$\{ x\sq\om: x\cap [k_{n_i},k_{m_{i+1}})=\emptyset\}
\sq U_{n_i}$$  
But this means that
$(U_{n_i}:i\in C)$ is a $\gamma$-cover of $\{x: x\sq^*x_\delta\}$.
But the construction of $(\ww_n:n\in\om)$ guarantees that
we can define them on all $n$ so that 
$(U_n\in\ww_n: n<\om)$ is an $\gamma$-cover of 
$$\finsub\cup\{x_\alpha:\alpha<\delta\} \cup \{x: x\sq^*x_\delta\}$$
which includes $X$.

\qed

\begin{cor}
In the iterated superperfect model we have
\par \centerline{(strong $\gamma$)- BC and not( $\gamma$-BC)}
\end{cor}


\begin{thebibliography}{99}


\bibitem{bd}
Baumgartner, James E.; Dordal, Peter; 
Adjoining dominating functions. 
J. Symbolic Logic 50 (1985), no. 1, 94--101.

\bibitem{dhm}
Dzamonja, M.; Hrusak, M.; Moore, J.;
Parameterized $\Diamond$ principles, eprint Feb 2003.

\bibitem{gm}
Galvin, Fred; Miller, Arnold W. $\gamma$-sets and other singular sets of real
numbers.  Topology Appl. 17 (1984), no. 2, 145--155.

\bibitem{gn} 
Gerlits, J.; Nagy, Zs.; Some properties of $C(X)$. I. Topology
Appl. 14 (1982), no. 2, 151--161. 

\bibitem{hechler}
Hechler, Stephen H.;
On the existence of certain cofinal subsets of $\sp{\om }\om $. 
Axiomatic set theory (Proc. Sympos. Pure Math., Vol. XIII, Part II, Univ.
California, Los Angeles, Calif., 1967), pp. 155--173.  Amer. Math. Soc.,
Providence, R.I., 1974. 

\bibitem{jmss}
Just, Winfried; Miller, Arnold W.; Scheepers, Marion; Szeptycki, Paul J.; The
combinatorics of open covers. II. Topology Appl. 73 (1996), no. 3, 241--266. 

\bibitem{laver}
Laver, Richard; 
On the consistency of Borel's conjecture. 
Acta Math. 137 (1976), no. 3-4, 151--169.

\bibitem{super}
Miller, Arnold W.; Rational perfect set forcing. Axiomatic set theory (Boulder,
Colo., 1983), 143--159, Contemp. Math., 31, Amer. Math. Soc., Providence, RI,
1984. 

\bibitem{truss}
Truss, John;
Sets having calibre $\aleph \sb{1}$. 
Logic Colloquium 76 (Oxford, 1976), pp. 595--612. 
Studies in Logic and Found. Math., Vol. 87, 
North-Holland, Amsterdam, 1977. 

\bibitem{tsaban}
Tsaban, Boaz; 
Strong gamma-sets and other singular spaces, eprint
arxiv.org math.LO/0208057.

\end{thebibliography}
\end{document}